\documentclass[a4paper]{amsart}
\usepackage{amssymb}
\input xy
\xyoption{all}

\newcommand{\R}{\mathbb R}

\newcommand{\N}{\mathbb N}
\newcommand{\T}{\mathcal T}

\newtheorem{theorem}{Theorem}

\newtheorem{problem}{Problem}

\begin{document}

\title{On monotonicity of comonotonically maxitive functional}


\author{Taras Radul}

\maketitle

Institute of Mathematics, Kazimierz Wielki University in Bydgoszcz, Poland;
\newline
Department of Mechanics and Mathematics, Ivan Franko National University of Lviv,
Universytettska st., 1. 79000 Lviv, Ukraine.
\newline
e-mail: tarasradul@yahoo.co.uk

\textbf{Key words and phrases:}  Fuzzy integral, comonotonically maxitive functional, monotone functional

\subjclass[MSC 2020]{ 28E10}

\begin{abstract} The comonotonic maxitivity property of functionals frequently appears in the characterization of fuzzy integrals based on the maximum operation. In some special cases, comonotonic maxitivity implies monotonicity of functionals. The question of whether this implication holds in general was posed by T. Radul (2023). It was shown in that paper that the implication is valid for finite compacta. In this article, we provide a negative answer to the general problem and discuss additional properties that need to be imposed to ensure the implication holds.
\end{abstract}

\maketitle

\section{Introduction}

In fact, most applications of non-additive measures in game theory, decision theory, economics, and related fields do not focus on measures as set functions, but rather on integrals that allow for the computation of expected utility or expected payoff. Several types of integrals with respect to non-additive measures have been developed for different purposes (see, for example, \cite{Grab}, \cite{KM}, \cite{LMOS}, and \cite{Den}). These integrals are commonly referred to as fuzzy integrals. The most well-known among them are the Choquet integral, which is based on addition and multiplication operations \cite{Ch}, and the Sugeno integral, which relies on maximum and minimum operations \cite{Su}. If the minimum operation is replaced by a general t-norm, the result is a generalization of the Sugeno integral known as a t-normed integral \cite{Sua}.

One of the central problems in the theory of fuzzy integrals is the characterization of such integrals as functionals on suitable function spaces (see, for example, Subchapter 4.8 in \cite{Grab}, which is devoted to characterizations of the Choquet and Sugeno integrals). A characterization of t-normed integrals was obtained in \cite{CLM} for finite compacta, and extended to the general case in \cite{Rad}.

It was observed in \cite{Grab} that in the special case of the Sugeno integral over finite sets, some of the conditions used in the characterization theorem of \cite{CLM} - notably monotonicity - are redundant, and a simpler characterization was proposed. This result was later generalized to the case of t-normed integrals on compacta in \cite{Rad1}. In particular, it was shown that monotonicity follows from other properties of the t-normed integral, among them the comonotonic maxitivity property.

The question of whether comonotonic maxitivity alone implies monotonicity was explicitly posed in \cite{Rad1}, where a positive answer was established for finite compacta. In the present paper, however, we demonstrate that the implication does not hold in general.

\section {Main result}

In what follows, all spaces are assumed to be compacta (compact Hausdorff space) except for $\R$ and the spaces of continuous functions on a compactum. All maps are assumed to be continuous. We shall denote the
Banach space of continuous functions on a compactum  $X$ endowed with the sup-norm by $C(X)$. For any $c\in\R$ we shall denote the
constant function on $X$ taking the value $c$ by $c_X$. We also consider the natural lattice operations $\vee$ and $\wedge$ ( on $C(X)$ and  its sublattice  $C(X,[0,1])$.

Let $X$ be a compactum.  We call two functions $\varphi$, $\psi\in C(X)$ comonotone (or equiordered) if $(\varphi(x_1)-\varphi(x_2))\cdot(\psi(x_1)-\psi(x_2))\ge 0$ for each $x_1$, $x_2\in X$. Let us remark that a constant function is comonotone to any function $\psi\in C(X)$.

Following \cite{Grab} we say that a functional  $I:C(X,[0,1])\to[0,1]$ is comonotonically maxitive if it preserves  $\vee$ for  comonotone functions.  It is known  \cite{Rad1} that monotonicity does not imply comonotonic maxitivity. It is easy to see that a comonotonically maxitive functional is monotone for  pairs of comonotone functions.

\begin{theorem}\cite{Rad1}\label{finite1} Let $X$ be a finite set and a  functional $\mu:C(X,[0,1])=[0,1]^n\to[0,1]$ be comonotonically maxitive. Then $\mu$ is monotone.
\end{theorem}

However, the following example shows that this implication does not hold in the general case.

Consider the compactum $X=\{p\}\cup A\subset\R$, where $A=\{1-\frac{1}{n}|n\in\N\}\cup\{1\}$ and $p=2$. For $t\in\{0,1\}$ define the function $f_t\in C(X,[0,1])$ by setting $f_t(p)=t$ and $f_t(a)=a$ for each $a\in A$. Define the following subsets of $C(X,[0,1])$: 

\begin{itemize}
\item $F_1=\{f\in C(X,[0,1]|f\le f_1\}$;
\item $F_2=\{f\in C(X,[0,1])|\max_{x\in X}f(x)\le f(p)\}$;
\item $F_3=\{f\in C(X,[0,1])|\max_{x\in X}f(x)<1\}$.
\end{itemize}  Put $G=F_1\cap (F_2\cup F_3)$. Define the functional $\nu:C(X,[0,1])\to\R$ by setting $\nu(f)=0$ if $f\in G$ and  $\nu(f)=1$ otherwise.

\begin{theorem}\label{general} The functional $\nu$ is comonotonically maxitive but not monotone.
\end{theorem}

\begin{proof} We clearly have  $f_0\le f_1$, but $\nu(f_1)=0<1=\nu(f_0)$. Hence, $\nu$ is not monotone.

Now consider any comonotone functions $f,g\in C(X,[0,1])$. It suffices to show that $\nu(f\vee g)=0$ whenever $f,g\in G$, and $\nu(f\vee g)=1$ otherwise.

Assume $f$, $g\in G$. Evidently,  $f\vee g\in F_1$. 

First, suppose $f$, $g\in F_2$. It is easy to see that $f\vee g\in F_2$.

Next, suppose  $f\in F_2$ and $g\in F_3$. Put $m=\max_{x\in X}g(x)<1\}$. 
\begin{itemize}
\item If $f(p)\ge m$, then $f\vee g\in F_2$.
\item If $f(p)\le m$, then $f\vee g\in F_3$.
\end{itemize}

Finally, if $f$, $g\in F_3$, then $f\vee g\in F_3$ for any $f$, $g\in F_3$. 

In all the above cases  $f\vee g\in G$ and hence $\nu(f\vee g)=0$.

Now suppose $f\notin G$. If $f\notin F_1$, then $f\vee g\notin F_1\supset G$. Consider the case  $f\in F_1\setminus (F_2\cup F_3)$. Then $f(p)<f(1)=1$. If $g\notin F_1$, then again $f\vee g\notin F_1\supset G$. Now assume  $g\in F_1$. If $g(p)=1$, take $a\in A\setminus\{1\}$ such that $f(a)>f(p)$. But since $g(a)\le f_1(a)=a<1=g(p)$ and this contradicts  the assumption that $f$ and $g$ are comonotone. Thus, $g(p)<1$ and we have $$(f\vee g)(p)<1=(f\vee g)(1),$$ which implies $f\vee g\notin F_2\cup F_3\supset G$. Therefore in all cases where  $f\notin G$ we have $f\vee g\notin G$, hence  $\nu(f\vee g)=1$.
\end{proof}

Let us note that if we impose some additional conditions on the functional
$\nu$, the implication becomes valid. 
Recall that a triangular norm $\ast$ is a binary operation on the closed unit interval $[0,1]$ which is associative, commutative, monotone and $s\ast 1=s$ for each  $s\in [0,1]$ \cite{PRP}.   We consider only continuous t-norms in this paper. Let $\ast$ be a continuous t-norm. For a compactum $X$, we denote  by $\T^\ast(X)$ the set of functionals $\mu:C(X,[0,1])\to[0,1]$ satisfying the following conditions:

\begin{enumerate}
\item $\mu(c_X)=c$ for each $c\in[0,1]$ ($\mu$ is normalized);
\item $\mu$ is comonotonically maxitive;
\item $\mu(c_X\ast\varphi)=c\ast\mu(\varphi)$ for each $c\in[0,1]$ and $\varphi\in C(X,[0,1])$ ($\mu$ is $\ast$-homogeneous).
\end{enumerate}

It follows from the results of \cite{Rad1} that each functional in $\T^\ast(X)$ is monotone.

Since the previously constructed functional  $\nu$  is not normalized, we may pose  the following intermediate problem:

\begin{problem} Let $X$ be a compactum and  let $\mu:C(X,[0,1])\to[0,1]$ be a  functional that is comonotonically maxitive and normalized. Is $\mu$ necessarily  monotone?
\end{problem}

\end{document}